\documentclass[10pt]{article}
\usepackage{graphicx}
\usepackage{todonotes}

\def\R{{\rm I\! R}}

\def \pTW*{\partial_{W } T}
\def \pTnW*{\partial^{(n)}_{W } T}

\newtheorem{theorem}{Theorem}
\newtheorem{corollary}{Corollary}

\title{On the  fixed points set of differential systems reversibilities} 

\author{Marco Sabatini 
\footnote{Dip. di Matematica, Univ. di Trento, I-38123 Povo, (TN) - Italy.
Phone: ++39(0461)281670, Fax: ++39(0461)281624, Email: marco.sabatini@unitn.it - \ \ \ \ \ \ \ \  \ \ \ \ \ \ \ \ \ \ \ \ This paper was partially supported by the GNAMPA, {\it Gruppo Nazionale per l'Analisi Matematica, la Probabilit\`a e le loro Applicazioni}. }
}
\date{ October 5, 2015}
\begin{document}
\maketitle
\begin{abstract}
We extend a result proved in \cite{Col} for mirror symmetries of planar systems to measure-preserving non-linear reversibilities of $n$-dimensional systems, dropping the analyticity and nondegeneracy conditions.  
\end{abstract}



\section{Introduction}

Let us consider a planar differential system
\begin{equation} \label{sys}
\dot z = F(z),
\end{equation}  
where $F(z)=(F_1(z),F_2(z)) \in C^1 (\Omega,\R^2)$,  $\Omega \subset \R^2$ open and connected, $z \in \Omega$. We denote by $\phi(t,z)$ the local flow defined by (\ref{sys}) in $\Omega$.   We assume the origin $O$ to be  a critical point of (\ref{sys}). One of the classical problems considered in the study of such systems is the so-called {\it center-focus} problem for a rotation point. It consists in  distinguishing among the two following possibilities,
\begin{itemize}
\item $O$ is an attractor;
\item $O$ is an accumulation point of cycles, in particular a center.
\end{itemize}
Poincar\'e and Lyapunov developped a procedure that applies to analytic systems with a non-degenerate critical point $O$, allowing to discern among the only possible cases for analytic systems, i. e. a center or a focus.

 If the system is non-analytic or the critical point is degenerate, in order to prove that a system has a center at $O$ one can either look for a first integral with an isolated extremum at $O$, or apply a symmetry argument first used by Poincar\'e, looking for a mirror symmetry of the solution set. The symmetry method is easy to apply if one knows the  symmetry line, since after a suitable axes rotation the symmetry conditions reduce to a parity verification on the components of the vector field. On the other hand, if such a line is unknown the question is obviously more complex, depending on a free parameter \cite{BL,Con}. In \cite{Col}, assuming the system to be analytic and the critical point $O$ to be non-degenerate, it was proved that  if a symmetry line exists, then it is contained in the zero-divergence set. This allows for a fast analysis of the possible mirror symmetries, since they are determined by the lowest order non-linear  terms of $F(z)$. In \cite{Col} this was applied to study quadratic and cubic polynomial systems, also considered in \cite{vW}.

The existence of a center  can also be proved finding a non-linear symmetry. We say that a diffeomorphism $\sigma \in C^1(\Omega,\Omega)$ of $\Omega$ onto itself is an {\it involution} if
$$
\sigma^2(z) = \sigma(\sigma(z)) = z \qquad  \forall z \in \Omega.
$$
An involution $\sigma$  is said to be a {\it reversibility} of the system (\ref{sys}) if
$$
\phi(t,\sigma(z)) = \sigma(\phi(-t,z)),
$$
 for all $t\in \R$ and $z\in\Omega$ such that both sides are defined. The non-linear generalization of the symmetry line of a mirror symmetry is the set $\delta$ consisting of the $\sigma$-fixed points, $\{ z \in \Omega : z = \sigma(z)\}$. Poincar\'e's argument can be easily extended to reversibilities, proving that an orbit intersecting $\delta$ is a cycle. As a consequence, if $\delta$ is a curve containing a rotation point of a reversible system, then such a point is a center. 
Reversibilities can be revealed by examining the symmetry properties of the vector field. In fact,  if $\sigma \in C^1(\Omega,\Omega)$, then the $\sigma$-reversibility of (\ref{sys}) can be checked by verifying the relationship
\begin{equation} \label{Jrev}
 V(\sigma(z)) = - J_\sigma(z) \cdot V(z) ,
\end{equation} 
where $J_\sigma(z)$ is the Jacobian matrix of $\sigma$ at $z$.  Finding a non-linear reversibility is obviously much more difficult than proving mirror symmetry with respect to a line. Moreover the fixed-points set $\delta$ of a non-linear reversibility $\sigma$ is not, in general, a line.
The existence of non-linear reversibilities in a neighbourhood of a nilpotent critical point was studied in \cite{T,BM,SZ}. Other results on reversibility, or using reversibility can be found in \cite{ACGG,CS,CRZ,GM,GV,L,LH,R,Z}.
 
 \smallskip
 
We say that $\sigma$ is {\it measure-preserving} if
 $$
 \mu(\sigma(A)) = \mu(A),
 $$
 where $\mu$ is the Lebesgue measure and $A \subset \Omega$ a measurable set. Measure-preserving diffeomorphisms appear frequently in different areas of mathematics. For instance, the flow defined by a Hamiltonian system is measure-preserving. Also, the maps involved in the celebrated Jacobian Conjecture are area-preserving \cite{BCW}. In fact, they are canonical transformations of the plane onto itself. 
The definitions of involution, reversibility and measure-preserving involution can be given in the same terms in dimension $n$. The  characterization (\ref{Jrev}) holds as well for $n$-dimensional systems.

 \smallskip
 
In this paper we prove that the fixed-points set of a  measure-preserving reversibility of a $n$-dimensional system is contained in the zero-divergence set of the vector field. We neither assume the system to be analytic, nor a critical point to be non-degenerate. In fact, our proof holds even for systems without critical points, hence it is also applicable to prove the existence of period annuli not contained in central regions. We apply our result to give an integrability condition for a class of planar $y$-quadratic systems. \\

 \bigskip
 
\section{The result}

Without changing notation, in this section we denote again by
\begin{equation} \label{sysn}
\dot z = F(z),
\end{equation}  
a differential system in an open and connected set $\Omega \subset \R^n$, where $F(z)=(F_1(z),\dots, F_n(z)) \in C^1 (\Omega,\R^n)$, $z=(z_1, \dots, z_n)$. Similarly, we write  $\phi(t,z)$ for the local flow defined by (\ref{sys}) in $\Omega$.   
 Let us denote by ${\rm div\,} F(z)$ the divergence of the vector field $F(z)$. We set
 $$
 \Delta_0 = \{ z \in \Omega : {\rm div\,} F(z) = 0\}.
 $$
 In  \cite{S} it was proved, for planar $\sigma$-reversible systems, that  every arc of orbit  containing both $z$ and $\sigma(z)$, contains a $\sigma$-fixed point, too. Such a proof extends  without changes to $n$-dimensional systems. In next theorem we show that  if $\sigma$ is a measure-preserving reversibility, then every  arc  of orbit  containing both $z$ and $\sigma(z)$, contains a zero-divergence point, too. As a consequence, we show that a regular $\sigma$-fixed point is a zero-divergence point. 
  
\begin{theorem} \label{divfix} Let $\sigma \in C^1(\Omega,\Omega)$ be a measure-preserving reversibility of (\ref{sys}). Then
\begin{itemize}
\item[i)] if there exist $z \in \Omega$ and $t_z  \in \R$ such that $\sigma(z ) = \phi(t_z ,z ) \neq z$, then the arc of orbit $\gamma_z$ connecting $z $ to $\sigma(z )$ contains a point $z_0 \in \Delta_0$.;
\item[ii)] if $z$ is a regular $\sigma$-fixed point, then $z \in \Delta_0$.
\end{itemize}
\end{theorem}
{\it Proof.}  
i) Possibly exchanging $z$ and $\sigma(z)$, we may assume  $t_z  >0$. 

By absurd, assume  the divergence  not to vanish on the arc of orbit $\gamma_z $ connecting  $z $ to $\sigma(z )$.  Then ${\rm div\,} F(z)$ has constant sign, say ${\rm div\,} F(z) >0$, and by the compactness of such an arc, $\gamma_z $ has an open neighbourhood $U $ such that ${\rm div\,} F(z) >0$ on $U $. By the continuity of $\phi(t,\cdot)$ there exists a measureable neighbourhood $A $ of $z $ such that $\mu(A ) > 0$ and $\phi(t,A ) \subset U $ for all $t \in [0,t _z]$. 
By Liouville theorem, in a positive divergence region the measure increases along the  local flow, hence
$$
t_1 , t_2 \in [0,t _z], \quad  t_1 < t_2 \ \  \Rightarrow \ \ \mu(\phi(t_1,A )) < \mu(\phi(t_2,A )) . 
$$
 On the other hand, for $ t \in (0,t _z] $ one has
 $$
 \mu(A )   <  \mu(\phi(t,A)) = \mu(\sigma(\phi(t,A))) =  \mu(\phi(-t,\sigma(A))) \leq  \mu(\sigma(A)) = \mu(A),
 $$
 contradiction.
 \\

ii) Since $F(z) \neq 0$, there exists $\varepsilon > 0$ such that $\phi(t,z)$ is a simple (i. e. injective) regular curve defined in the interval $[-\varepsilon,\varepsilon]$. Then
$$
\sigma(\phi(-\varepsilon,z)) = \phi(\varepsilon,\sigma(z)) = \phi(\varepsilon,z),
$$
The hypotheses of point $i)$  hold, since the orbit through $z$ contains the distinct points $\phi(-\varepsilon,z) $ and $\sigma(\phi(-\varepsilon,z)) = \phi(\varepsilon,z)$, hence such an arc contains a point of $\Delta_0$. Similarly, every arc $\phi([-\eta,\eta],z)$, with $0 < \eta < \varepsilon$, contains a point of $\Delta_0$. By the continuity of ${\rm div \, } F$, this implies that ${\rm div \, } F(z) = 0$.
\hfill  $\clubsuit$  \\

If the orbit $\gamma_z$ is a cycle, then both its arc going from $z$ to $\sigma(z) \neq z$, and the second arc going from $\sigma(z)$  to $z$ contain a point of $\Delta_0$. In \cite{S} it was proved that every every cycle $\gamma$ of a $\sigma$-reversible system satisfying $\sigma(\gamma) = \gamma$ contains exactly two $\sigma$-fixed points.

\bigskip

In general the converse to the statement of theorem \ref{divfix} is not true, there exist reversibilities with zero-divergence points  that are not fixed points of the reversibility. For instance, the system
$$\left\{
\matrix{\dot x = y  \hfill \cr \dot y = -x\hfill}
\right.
$$
is hamiltonian, hence has identically zero divergence, and most of its regular points are not fixed points of the reversibility $\sigma(x,y)=(-x,y)$.  \\

Next corollary's proof is immediate.

\begin{corollary} \label{mirror} Let $\sigma \in C^1(\Omega,\Omega)$ be a mirror symmetry of (\ref{sys}). 
Then the symmetry hyperplane is contained in $\Delta_0$.
\end{corollary}

\section{An integrability condition for some $y$-quadratic systems}

In this section we consider a class of planar $y$-quadratic systems,
\begin{equation}   \label{yquadra}
\left\{
\matrix{\dot x = r(x)y + s(x) \hfill \cr \dot y = -g(x) - f(x)y - h(x)y^2 .\hfill}
\right.
\end{equation}
If $r(x) \neq 0$, such a system is orbitally equivalent to a second order scalar O.D.E. Its zero-divergence set has a very simple form, since it is the graph of a one-variable function, $y=\alpha(x)$. 
Since a  class of area-preserving  involutions is provided by the family of triangular maps, 
\begin{equation}   \label{sa}
\sigma(x,y) = (x, 2\alpha(x) - y),
\end{equation}
having the curve $y = \alpha(x)$ as fixed curve, we look for conditions under which (\ref{sa}) is a reversibility for (\ref{yquadra}). 
This suggests to perform a change of variables, taking $y=\alpha(x)$ into one axis, and check the reversibility conditions in the new variables. It comes out that, under the hypotheses of next theorem,  the new conditions are equivalent to the classical mirror symmetry of a vector field.

\begin{theorem} \label{teoquadra}Let $r, s, f, g, h \in C^1(a,b)$, with $a < 0 < b$. If  for all $x \in (a,b)$, 
$$
r(x) \neq 0, \qquad r'(x) - 2h(x) \neq 0,
$$
and
\begin{equation}   \label{cond}
r(x)f(x) - 2h(x)s(x) - s'(x)r(x) + s(x)r'(x) = 0,
\end{equation}
then the system (\ref{yquadra}) is $\sigma$-reversible, with 
$$
\sigma(x,y) =    \left(x, 2\,  \frac { f(x) - s'(x)} {r'(x) - 2h(x)} - y \right) = \left(x, -  \frac {2 s(x)} {r(x)} - y \right).
$$
\end{theorem}
{\it Proof.}   Setting $F(x,y) = ( r(x)y + s(x) , -g(x) - f(x)y - h(x)y^2)$, one has
$$
 {\rm div\,} F(x,y) = r'(x)y + s'(x)  - f(x) - 2h(x)y. 
$$
One has  $ {\rm div\,} F(x,y) = 0$ on the curve $y=\alpha(x)$, where
$$
\alpha(x) = \frac {f(x) - s'(x)} {r'(x) - 2h(x)}.
$$

The change of variables $(u,v)=\Lambda(x,y)=(x,y-\alpha(x))$ takes the system (\ref{yquadra}) into the system
\begin{equation}   \label{yquadrauv}
\left\{
\matrix{\dot u = r v + r  \alpha  + s \hfill \cr \dot v = - g -  f \alpha -h \alpha^2 - \alpha \alpha' r -\alpha' s - (f + 2 h \alpha + \alpha' r) v - h v^2 , \hfill}
\right.
\end{equation}
omitting the dependence on $u$. The curve $y=\alpha(x)$ is taken by $\Lambda$  into the $u$ axis.

\smallskip

The flow $\phi$ defined by (\ref{yquadra}) is $\sigma$-reversible if and only if the flow $\phi_\Lambda$ defined by (\ref{yquadrauv}) is reversible with respect to the involution
$$
\sigma_\Lambda(u,v) = \Lambda (\sigma(\Lambda^{-1}(u,v))) = (u,-v).
$$
Hence the system (\ref{yquadra}) is $\sigma$-reversible if and only if the components of the system (\ref{yquadrauv}) satisfy the usual parity properties for  mirror reversibility with respect to the $u$ axis,
$$
\dot u(u,-v) = -\dot u(u,v), \qquad  \dot v(u,-v) = \dot v(u,v) .
$$
This is equivalent to the conditions
\begin{equation}   \label{parcond}
\left\{
\matrix{s+ r \alpha   = 0\hfill \cr   f + 2h \alpha +  r \alpha' = 0   . \hfill}
\right.
\end{equation}
The condition (\ref{cond}) can be written as
$$
s = \frac {r(f-s')} {2h - r'},
$$
hence
$$
s + r \alpha  = \frac {r(f-s')} {2h - r'} +  \frac { r(f - s')} {r' - 2h} = 0,
$$
which shows that the first equation in (\ref{parcond}) is satisfied. As for the second one, from the first equation at the points where $r(x) \neq 0$ one has
$$
\alpha = - \frac sr.
$$
Hence one can write the second equation as follows
$$
f + 2h \alpha +  r \alpha' = f - \frac {2hs}{r} - r \frac {s'r - r's }{r^2} = \frac {fr - 2hs - s'r + r's}{r^2} = 0,
$$
because the numerator vanishes by (\ref{cond}). 
\hfill $\clubsuit$ \\

As a consequence of the above theorem, on has the following application to the existence of centers. The proof is immediate.

\begin{corollary} \label{mirror} Under the hypotheses ot theorem \ref{teoquadra}, if the origin $O$ is an isolated critical point of rotation type, then $O$ is a center.
\end{corollary}

The inequalities in the hypotheses of theorem \ref{teoquadra} are not very restrictive. In fact, from the first equation of (\ref{parcond}) follows that under the hypotheses of theorem \ref{teoquadra},  $r(x_0) = 0$ implies $s(x_0) = 0$. 
In this case the line $x = x_0$ is an invariant line for the system (\ref{yquadra}), hence its points do not belong to a central region. 

\smallskip

Moreover, if $r'(x_0) - 2h(x_0) = 0$ and $f(x_0) - s'(x_0) \neq 0$, then the curve $y = \alpha(x)$ has a vertical asymptote at $x_0$, hence the transformation $\Lambda$ cannot be extended continuously to a transformation defined also at $x_0$. 

If $r'(x) - 2h(x) \equiv 0$ on an interval $I$, then the system divergence in the set $I \times \R$ depends only on $x$, hence $\Delta_0$ is a family of lines parallel to the $y$ axis. Since $\Delta_0$ contains the reversibility curve, this is not compatible with a triangular reversibility of the form $\sigma(x,y) = (x,2\alpha(x) - y)$.

\smallskip

Finally, we observe that under the hypotheses of theorem \ref{teoquadra} the reversibility curve is the vertical isocline.
\smallskip

\smallskip

We give two examples of $\sigma$-reversible systems equivalent to second order differential equations. Both satisfy the hypotheses of theorem  \ref{teoquadra} and have a center at the origin. 

\smallskip

\begin{figure}[ht!] 
\caption{Some orbits of system (\ref{esuno})}  \label{fig1}
\centering
	\includegraphics[scale=0.5]{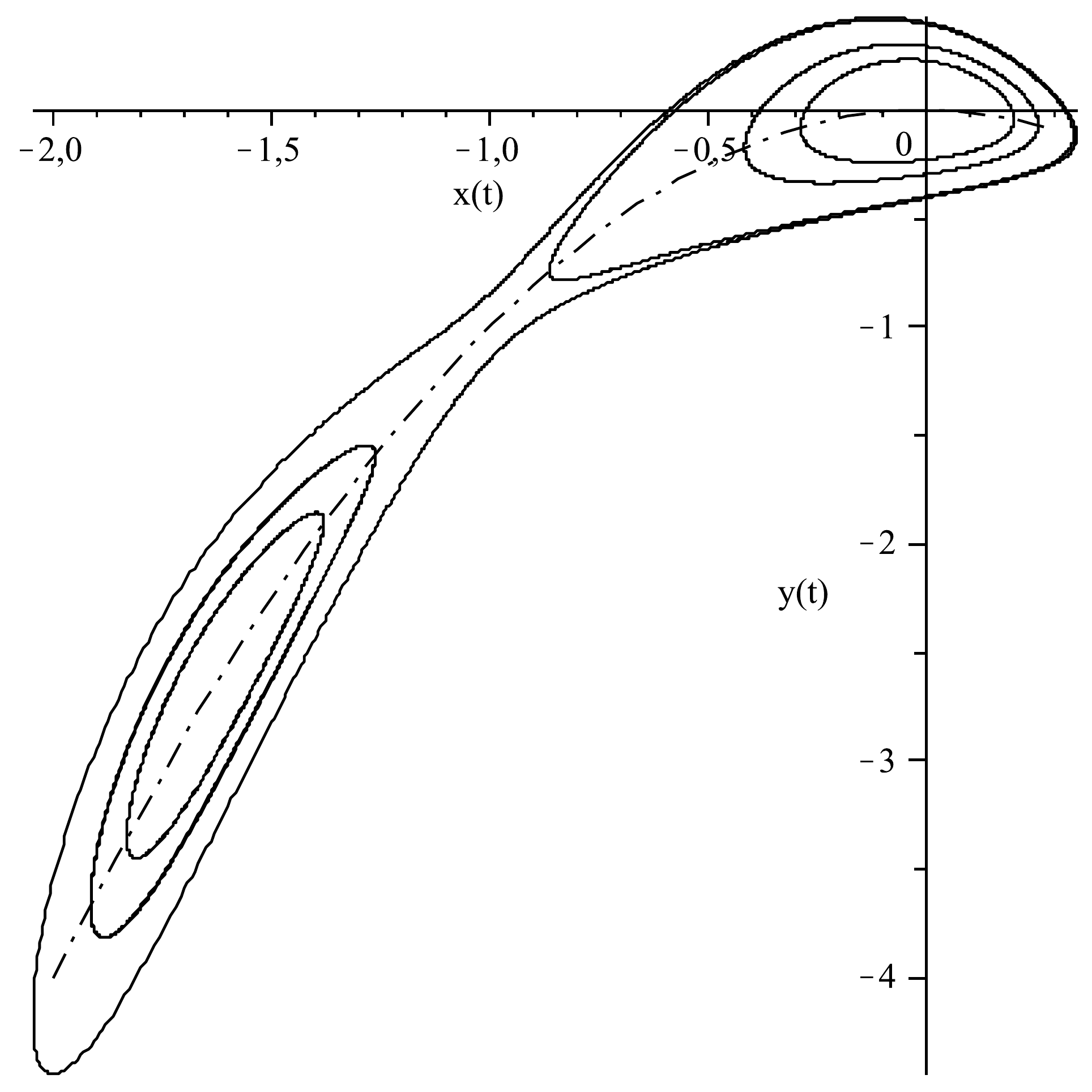}
\end{figure}

In this example we choose $r(x) \equiv 1$, $s(x) = x^2$, $g(x) = x$, $f(x) = 2x + 2x^2$, $h(x) \equiv 1$. In this case one has
$\alpha(x) = -x^2$.
\begin{equation}   \label{esuno}
\left\{
\matrix{\dot x = y +  x^2 \hfill \cr \dot y = -x - 2xy(1+x) - y^2 . \hfill}
\right.
\end{equation}
 The local phase portrait is sketched in figure \ref{fig1}, where the dash-dotted curve is the $\sigma$-fixed point curve. 
\vskip0.5cm

\smallskip

\begin{figure}[ht!] 
\caption{Some orbits of system (\ref{esdue})}  \label{fig2}
\centering
	\includegraphics[scale=0.5]{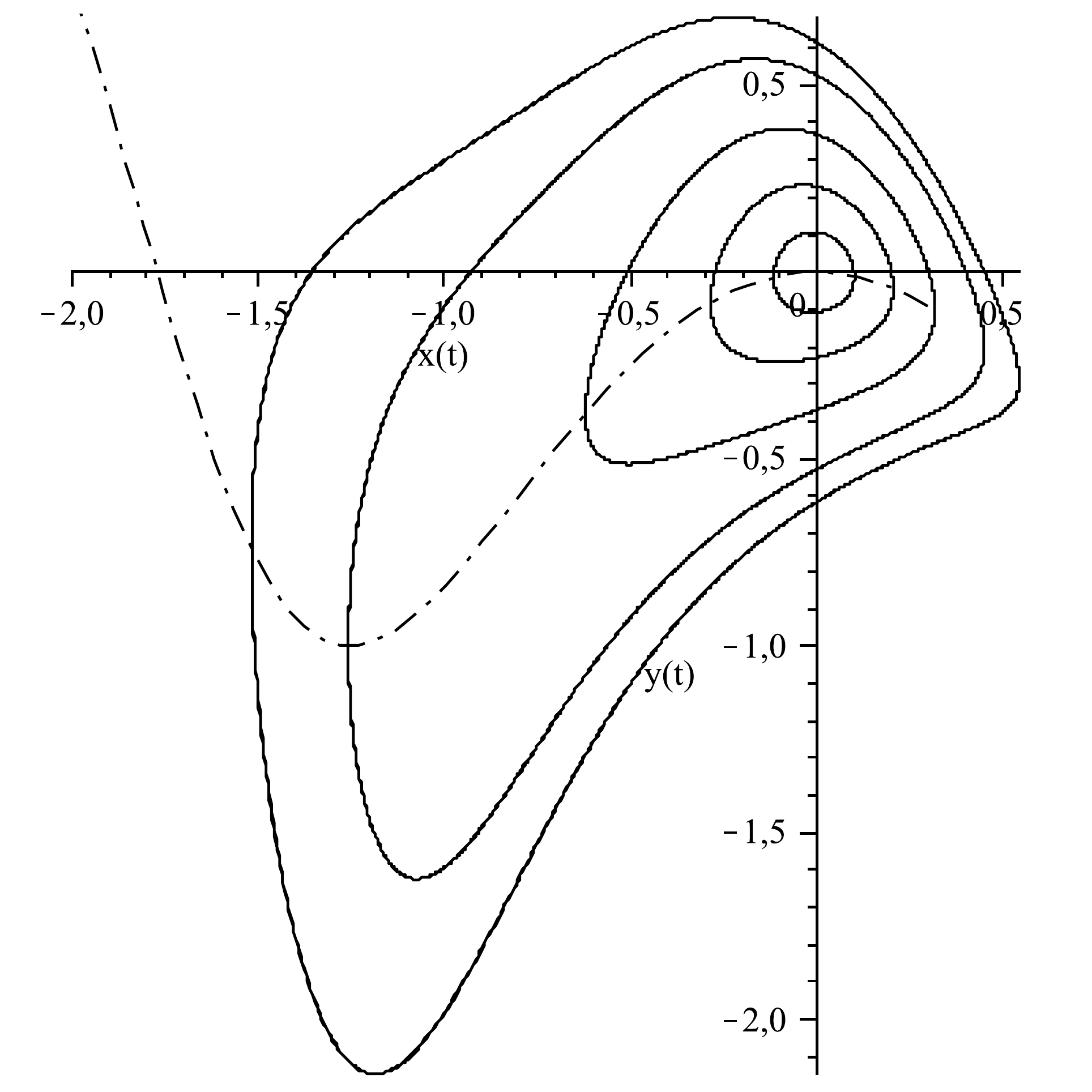}
\end{figure}

In next example we choose $r(x) \equiv 1$, $s(x) = \sin x^2$, $g(x) = x$, $f(x) = 2\sin x^2 + 2x \cos x^2$, $h(x) \equiv 1$. In this case one has $\alpha(x) = -\sin(x^2)$.
\begin{equation}   \label{esdue}
\left\{
\matrix{\dot x = y +  \sin x^2 \hfill \cr \dot y = - x - (2 \sin x^2 + 2x \cos x^2)y - y^2 . \hfill}
\right.
\end{equation}
 The local phase portrait is sketched in figure \ref{fig2}, where the dash-dotted curve is the $\sigma$-fixed point curve. 


\vskip0.5cm



\begin{thebibliography}{99}


\bibitem {ACGG} {\sc A. Algaba, I. Checa, C. Garcia, E. Gamero}, {\it On orbital-reversibility for a class of planar dynamical systems},
Commun. Nonlinear Sci. Numer. Simul.   {\bf 20, 1}  (2015),  229 -- 239. 

\bibitem {BCW} {\sc H. Bass, E. H. Connell, D. Wright}, {\it The Jacobian conjecture: reduction of degree and formal expansion of the inverse},   {\bf 7, 2}  (1982),  287 -- 330. 

\bibitem {BM} {\sc M. Berthier, R. Moussu}, {\it R\'eversibilit\'e et classification des centres nilpotents}, Ann. Inst. Fourier (Grenoble)    {\bf 44, 2} (1994), 465 -- 494.

\bibitem {BL} {\sc T. R. Blows, N. G. Lloyd}, {\it The number of limit cycles of certain polynomial differential equations}, Proc. Roy. Soc. Edinburgh Sect. A,   {\bf 98} (1984),  215 -- 239.

\bibitem{CS} {\sc J. Chavarriga, M. Sabatini}, {\it A survey of  isochronous centers}, Qual. Theory Dyn. Syst.   {\bf 1, 1} (1999), 1 -- 70.  

\bibitem {CRZ} {\sc X. Chen, V. G. Romanovski, W. Zhang}, {\it Critical periods of perturbations of reversible rigidly isochronous centers},
 J. Differential Equations {\bf 251, 6} (2011),  1505 -- 1525.  
 
\bibitem {Col} {\sc C. B. Collins}, {\it Poincar\'e's reversibility condition}, J. Math. Anal. Appl.   {\bf 259, 1} (2001),  168 -- 187.

\bibitem{Con} {\sc R. Conti}, {\it Centers of planar polynomial systems. A review.}  Le Matematiche {\bf 53} (1998), 207 -- 240.  

\bibitem {GV} {\sc A. Garijo, J. Villadelprat}, {\it Algebraic and analytical tools for the study of the period function}, J. Differential Equations   {\bf 257, 7} (2014), 2464 -- 2484.

\bibitem {GM} {\sc J. Gin\'e, S. Maza}, {\it The reversibility and the center problem}, Nonlinear Anal.   {\bf 74, 2} (2011), 695 -- 704.

\bibitem {L} {\sc C. Liu}, {\it The cyclicity of period annuli of a class of quadratic reversible systems with two centers},  J. Differential Equations  {\bf 252, 10} (2012),  5260 -- 5273.

\bibitem {LH} {\sc C. Liu, M. Han}, {\it Bifurcation of critical periods from the reversible rigidly isochronous centers},  Nonlinear Anal.   {\bf 95} (2014),  388 -- 403.

\bibitem {R} {\sc V. G. Romanovski}, {\it Time-reversibility in 2-dim systems},  Open Syst. Inf. Dyn.    {\bf 15, 4} (2008),  359 -- 370.

\bibitem {S} {\sc M. Sabatini}, {\it Every period annulus is both reversible and symmetric},  Arxiv:1504.04530, available at http://arxiv.org/abs/1504.04530.
 
 nilpotent
 \bibitem {SZ} {\sc E. Strozina, H. Zoladek}, {\it The analytic and formal normal form for the nilpotent singularity},
J. Differential Equations {\it 179, 2} (2002), 479 -- 537.

nilpotent
\bibitem {T} {\sc F. Takens}, {\it Singularities of vector fields},
Inst. Hautes ƒtudes Sci. Publ. Math. {\it 43} (1974), 47 -- 100. 

\bibitem {vW} {\sc W. von Wahl}, {\it Remarks on lines of reversibility for Poincar\'e's centre problem},  Analysis   {\bf 29} (2009),  259 -- 264.

\bibitem {Z} {\sc H. Zoladek}, {\it The classification of reversible cubic systems with center},  Topol. Methods Nonlinear Anal.   {\bf 4, 1} (1994),  79 -- 136.




 

 \end{thebibliography}
\enddocument